\documentclass[a4paper,11pt,reqno]{amsart}
\usepackage{amsmath}
\usepackage{amssymb}
\usepackage{amsthm}
\usepackage{nicefrac}
\usepackage{nccmath}
\usepackage{enumerate}
\usepackage{fixltx2e}
\usepackage[utf8]{inputenc}
\usepackage[english]{babel}
\usepackage{geometry}
\usepackage{graphicx}
\usepackage{subcaption}
\usepackage{relsize}
\usepackage{yfonts}
\usepackage{calligra}
\usepackage{amscd}
\usepackage{amssymb, latexsym}
\usepackage{mathrsfs}
\usepackage{fancyhdr}
\usepackage{enumerate}
\usepackage{setspace}	
\usepackage{latexsym}

\numberwithin{equation}{section}
\theoremstyle{plain}
\newtheorem{thm}{Theorem}[section]
\newtheorem{lem}[thm]{Lemma}
\newtheorem{pro}[thm]{Proposition}

\newtheorem{defn}[thm]{Definition}

\newtheorem{exa}[thm]{Example}
\numberwithin{thm}{section}
\numberwithin{equation}{section}

\allowdisplaybreaks

\providecommand{\keywords}[1]
{
	\small	
	\textbf{\textit{Keywords---}} #1
}
\title{On Weyl-Heisenberg Frames}

\begin{document}
	\author[Satyapriya, R. Kumar, A. K. Sah and Sheetal]{Satyapriya$^1$, Raj Kumar$^2$, Ashok K. Sah$^3$ and Sheetal$^4$}
	\address[1]{Department of Mathematics, University of Delhi, Delhi, India.}
	\email{kmc.satyapriya@gmail.com}
	\address[2]{Department of Mathematics, Kirori Mal College, University of Delhi, Delhi, India.}
	\email{rajkmc@gmail.com}
	\address[3]{Department of Mathematics, University of Delhi, Delhi, India.}
	\email{ashokmaths2010@gmail.com}
	\address[4]{Department of Mathematics, SGTB Khalsa College,  University of Delhi, Delhi.}
	\email{sheetal@sgtbkhalsa.du.ac.in}
	\keywords{Frames, $K$-Frames, Weyl-Heisenberg Frames}
	\subjclass[2010]{Primary: 42C15, 46B15.}
	\thanks{$^1$Corresponding author. Email: kmc.satyapriya@gmail.com\\
	The work of the corresponding author is supported by the Senior Research Fellowship (SRF) of Human Resource Development Group, Council of Scientific and Industrial Research (HRDG-CSIR), India (Grant No: 09/045(1653)/2019-EMR-I)}
	\maketitle
	\begin{abstract}
		In this paper we have generalized and studied the $K$-Weyl-Heisenberg frames, where $K$ is a bounded linear operator	on $L^2(\mathbb{R}^d)$. We have obtained necessary and sufficient conditions for acertain system to be a $K$-Weyl-Heisenberg frame. We have also given the invariance property of these $K$-Weyl-Heisenberg frames.	
	\end{abstract}
	\section{Introduction}
		\label{intro}
		Today frame theory is a central tool in applied mathematics and
		engineering.  Frames are widely used in sampling
		theory, wavelet theory, wireless communication, signal processing, image processing ,
		differential equations, filter banks, wireless sensor
		network, many more. Frames for Hilbert spaces introduced by  Duffin and Schaeffer \cite{DS},
		to address difficult problems from the theory of nonharmonic Fourier series.
		The theory of frames re-emerged in 1986, when Daubechies, Grossmann and Meyer \cite{DGM} found new
		applications to wavelet and Gabor transforms in which frames played
		an important role.
		A system $\{f_k\}$ in a separable Hilbert space $\mathcal{H}$ with inner product $\langle ., . \rangle$ is
		called \emph{frame} (Hilbert) for $\mathcal{H}$ if there exists
		positive constants $A$ and $B$ such that
		\begin{align}\label{e1}
			A \|f\|^2\leq  \|\langle f, f_k\rangle\|^2_{\ell^2} \leq B \|f\|^2,
			\ \text{for \ all} \  f \in \mathcal{H}.
		\end{align}
		If upper inequality in \eqref{e1} holds, then we say that $\{f_k\}$ is
		a Bessel sequence for $\mathcal{H}$. The positive constants $A$ and
		$B$ are called \emph{lower} and \emph{upper frame bounds} of the
		frame $\{f_n\}$, respectively. They are not unique..
		An introduction to frames in applied mathematics and engineering can
		be found in a so nice book by Casazza and Kutynoik \cite{CK}. For basic theory
		in frames, an interested reader may be refer to \cite{C,OC,FS}.
		
		This paper is organized as follows: In Section $2$ we give basic
		definitions and results which will be used throughout this paper. In
		Section $3$, we first study the notion of  $K$-Weyl-Heisenberg frames (or $K$-$\mathcal{W} \ \mathcal{H}$ frames ) in $L^2(\mathbb{R}^d)$. Necessary and
		sufficient conditions for a certain system in $L^2(\mathbb{R}^d)$ to
		$K$-$\mathcal{W} \ \mathcal{H}$  frames for
		$L^2(\mathbb{R}^d)$ are obtained. Invariance property of
		$K$-Weyl-Heisenberg frames is given.
	\section{Preliminaries And Notations}
		Some basic notations are listed in this section. Throughout this paper, we use the following notations.\\
		For $1 \leq p < \infty$, let $L^p(\mathbb{R}^d)$ denote the
		Banach space of complex-valued Lebesgue integrable functions $f$ on
		$\mathbb{R}^d$ satisfying
		\begin{align*}
		\|f\|_p = \left( \int_{\mathbb{R}^d}|f(t)|^p dt \right)^\frac{1}{p} <
		\infty.
		\end{align*}
		For $p=2$, an  inner product on $L^p(\mathbb{R})$ is given by
		\begin{align*}
		\langle f, g \rangle = \int_{\mathbb{R}^d}f \overline{g}dt,
		\end{align*}
		where $\overline{g}$ denote the complex conjugate of $g$.
		
		We  define the unitary operators $T_a ,E_b; a, b\in \mathbb{R}^d$  on
		$L^p(\mathbb{R}^d)$ by :\\
		\textbf{Translation} $\leftrightarrow$   $T_af(t)=f(t-a)$.\\
		\textbf{Modulation} $\leftrightarrow$ $E_bf(t)=e^{2\pi ib.t}f(t).$\\
		For $g \in L^2(\mathbb{R})$, $E_{mb}T_{na}g(t)=e^{2\pi imbt}g(t-na)$. If $a, b > 0$ and $(g, a,
		b)= \{E_{mb}T_{na}g\}_{m, n \in \mathbb{Z}}$ is a frame for
		$L^2(\mathbb{R})$, then we call $(g, a, b)$ a \emph{Gabor frame} or
		a \emph{Weyl-Heisenberg frame} for $L^2(\mathbb{R})$. Casazza \cite{C} introduced and studied irregular
		Weyl-Heisenberg  frames for $L^2(\mathbb{R})$. Let
		$(x_m,y_n)\in \mathbb{R}^2$ and let $g \in L^2(\mathbb{R})$. A
		system of the form $\{E_{x_m}T_{y_n}g(t)\}_{m,n\in \mathbb{\mathbb{Z}}}$ is
		called an \emph{irregular Weyl-Heisenberg  frame} (or \emph{$IWH$
			frame}) for $L^2(\mathbb{R}),$  if $\{E_{x_m}T_{y_n}g(t)\}_{m,n\in
			\mathbb{\mathbb{Z}}}$ is a frame for $L^2(\mathbb{R})$.
		For further study one may refer to some recent work in Weyl-Heisenberg frames (modulation and translation) and its extension to wave packet decomposition (modulation, translation and dilation) in \cite{G,FS,KS,KS1,S,SV,SVR}
		
		\begin{defn}\cite{C1}
			A family of real numbers $\{\lambda_n\}_{n\in \mathbb{Z}}$ has
			\emph{uniform density} $D\equiv D(\{\lambda_n\})$ if there is an
			$L>0$ such that for all $n\in \mathbb{Z},$ we have
			$|\lambda_n-\frac{n}{D}|\leq L.$
		\end{defn}
		
		The following lemma provides necessary and sufficient conditions for
		a given system in a Hilbert space to be a Bessel sequence.
		
		\begin{lem}\cite{OC} Let $\{f_n\}_{n=1}^{\infty}$ be a sequence in
			$\mathcal{H}$  and $B>0$ be given. Then, $\{f_n\}_{n=1}^{\infty}$ is
			a Bessel sequence with Bessel bound $B$ if and only if \
			$T\{c_n\}_{n=1}^{\infty}=\sum_{n=1}^{\infty}c_nf_n$ defines a
			bounded linear operator from $\ell^2$ into $\mathcal{H}$ and
			$\|T\|\leq\sqrt{B}$.
		\end{lem}
		
		A class of bounded linear operator from a Banach space $\mathcal{X}$
		into a Banach space $\mathcal{Y}$ is denoted by
		$\mathcal{B}(\mathcal{X}, \mathcal{Y})$. If $\mathcal{X} =
		\mathcal{Y}$, then we write $\mathcal{B}(\mathcal{X}, \mathcal{X})
		\equiv \mathcal{B}(\mathcal{X}) $. Let $(\mathcal{H}, \langle .,
		.\rangle)$ be a Hilbert space. An operator $T \in
		\mathcal{B}(\mathcal{H})$ is said to be positive if $\langle Tf, f
		\rangle \geq 0$, for all $f \in \mathcal{H}$. In symbol we write $T
		\geq 0$. For $T_1, T_2 \in \mathcal{B}(\mathcal{H})$, if $T_1 - T_2 \geq 0$, then we write $T_1\geq T_2$.
		
		\begin{thm} \cite{RGD}\label{t2.3}
			Let $T_1\in \mathcal{B}(H_1,H), T_2\in \mathcal{B}(H_2,H)$ be two
			bounded linear operators, where $H, H_1, H_2$ are Hilbert spaces.
			The following statement are equivalent
			\begin{enumerate}[(i)]
				\item $R(T_1)\subset R(T_2)$, where $R(T_i)$ denote the range of
				$T_i$.
				\item $T_1T_1^{*}\leq \lambda ^2 T_2T_2^{*}$ for some $\lambda\geq0$
				\item There exists a bounded linear operator $S\in B(H_1,H_2)$ such
				that $T_1=T_2S .$
			\end{enumerate}
		\end{thm}

	\section{Weyl-Heisenberg Frames in $L^2(\mathbb{R}^d)$}
		\begin{defn}
			Let $g \in L^2(\mathbb{R}^d)$ , $\{C_m\}_{m\in \mathbb{Z}} \subset
			\mathbb{R}^d$, $ k\in {\mathbb{Z}}^d$ and $\mathbf{B}\in GL_d({\mathbb{R})}$, let $K \in \mathcal{B}(L^2(\mathbb{R}^d))$. A systen $\{E_{C_m}T_{\mathbf{B}n}g\}_{m\in \mathbb{Z}, n \in \mathbb{Z}^d} $ is said to be
			\emph{$K$-Weyl-Heisenberg frame} (or \emph{$K$-$\mathcal{W} \ \mathcal{H}$ frame}) for $L^2(\mathbb{R}^d)$
			if there exist positive constants
			$0 < \mathrm{A} \leq\mathrm{B}< \infty$ such that
			\begin{align}\label{e3.1}
				\mathrm{A}\|K^{*}f\|^2\leq \sum_{m\in \mathbb{Z}, n \in \mathbb{Z}^d}|\langle f,E_{C_m}T_{\mathbf{B}n}g\rangle|^2 \leq B\|K f\|^2, \  \text{for   all} \ f\in L^2(\mathbb{R}^d)
			\end{align}
		\end{defn}
		The positive constants $\mathrm{A}$ and $\mathrm{B}$ are called \textit{lower}
		and \emph{upper  frame bounds} of the $K$- Weyl-Heisenberg frame $\{E_{C_m}T_{\mathbf{B}n}g\}_{m\in \mathbb{Z}, n \in \mathbb{Z}^d} $ ,
		respectively . If upper ineqality  in \eqref{e3.1} is satisfied, then
		$\{E_{C_m}T_{\mathbf{B}n}g\}_{m\in \mathbb{Z}, n \in \mathbb{Z}^d} $ is called the
		\emph{$K$\textit{-Bessel sequence} for
			$L^2(\mathbb{R}^d)$} with Bessel bound $\mathrm{B}$.
		Suppose that $\mathcal{F} \equiv \{E_{C_m}T_{\mathbf{B}n}g\}_{m\in \mathbb{Z}, n \in \mathbb{Z}^d} $ is a $K$-$\mathcal{W} \ \mathcal{H}$ frame for $L^2(\mathbb{R}^d)$. The operator $V:l^2\bigoplus l^2\longrightarrow L^2(\mathbb{R}^d)$ given by
		\begin{align*}
			V\{c_{mn}\}_{m\in \mathbb{Z}, n \in \mathbb{Z}^d}=\sum_{m\in \mathbb{Z}, n \in \mathbb{Z}^d}c_{mn}E_{C_m}T_{Bn}g,
		\end{align*}
		is called the\textit{ pre-frame operator} or \textit{synthesis operator} associated with $\mathcal{F}$ and the adjoint operator
		$V^{*}:L^2(\mathbb{R}^d)\longrightarrow l^2\bigoplus l^2$ is given by
		\begin{align*}
			V^{*}\textit{f}=\{\langle f,E_{C_m}T_{\mathbf{B}n}g\rangle\}_{m\in \mathbb{Z}, n \in \mathbb{Z}^d}
		\end{align*}
		is called the \textit{analysis operator} associated with $\mathcal{F}$. Composing
		$V \ and \ V^{*}$, we obtain the \textit{frame operator}
		$\mathcal{S}:L^2(\mathbb{R}^d)\longrightarrow L^2(\mathbb{R}^d)$ given
		by
		\begin{align*}
			\mathcal{S}f=VV^{*}f=\sum_{m\in \mathbb{Z}, n \in \mathbb{Z}^d}\langle f,E_{C_m}T_{\mathbf{B}n}g\rangle E_{C_m}T_{Bn}g.
		\end{align*}
		Since $\mathcal{F}$ is a  $K$-Bessel sequence for
		$L^2(\mathbb{R}^d)$, the series defining $\mathcal{S}$ converges
		unconditionally for all $f\in L^2(\mathbb{R}^d).$  Note that, in
		general, frame operator of $\mathcal{F}$ is not invertible on
		$L^2(\mathbb{R}^d)$, but it is invertible on a subspace
		Range$(K)\subset L^2(\mathbb{R}^d)$. If Range$(K)$
		is closed , then the Pseudo-inverse $K^{\dag}$ of
		$K$ exist such that $KK^{\dag}f=f$, for all $f \in
		\text{Range}(K)$ and
		$KK^{\dag}|_{\text{Range}(K)}=I_{\text{Range}(K)}$,\\
		Therefore\\
		$\left(K^{\dag}|_{\text{Range}(K)}\right)^*K^{*}=I_{\text{Range}(K)}^{*}$.\\
		Now for any $f\in \text{Range}(K)$, we have,
		$\|f\|=\|\left(K^{\dag}|_{\text{Range}(K)}\right)^*K^{*}f\|\leq\|K^{\dag}\|\|K^{*}f\|.$\\
		Therefore
		\begin{align*}
			\langle \mathcal{S}f,f\rangle\geq A\|K^{*}f\|^2\geq
			A\|K^{\dag}\|^{-2}\|f\|^2,  \text{ for \ all } f \in
			\text{Range}(K).
		\end{align*}
		Again by the definition of a $K$-$\mathcal{W} \ \mathcal{H}$ frame, we get
		\begin{align*}
		A\|K^{\dag}\|^{-2}\|f\|^2\leq\|Sf\|\leq B \|f\|, \ \text{for \ all} f \in \text{Range}(K),
		\end{align*}
		Thus, the operator  $\mathcal{S}:\text{Range}(K)\longrightarrow
		\mathcal{S}(\text{Range}(K))$ is a homeomorphism.\\
		Moreover, we obtain
		\begin{align*}
			B^{-1}\|f\|\leq\|\mathcal{S}^{-1}f\|\leq
			A^{-1}\|K^{\dag}\|^2\|f\|, \  \text{for \ all }   f \in
			\mathcal{S}(\text{Range}(K)).
		\end{align*}
		
		\begin{exa}
			Let $\chi_{[0,1]}$ denote the characteristic  function for the
			interval [0,1], and let $C_m=m$, \ $\mathbf{B}=I$, where $m,n \in
			\mathbb{Z}$ and $g=\chi_{[0,1]}$. Choose $K=T_\xi$, a
			translation operator. Then, $\mathcal{F} \equiv\{e^{2\pi
				im(\bullet-\xi)}\chi_{[0,1]}(\bullet-\xi-n)\}_{m\in \mathbb{Z}, n \in \mathbb{Z}^d}
			$ is a $K$-$\mathcal{W} \ \mathcal{H}$-frame for
			$L^2(\mathbb{R}^d)$ with bounds $A = B = 1$.\\ Indeed
			\begin{align*}
			\sum_{m\in \mathbb{Z}, n \in \mathbb{Z}^d}|\langle
			f,e^{2\pi im(\bullet-\xi)}\chi_{[0,1]}(\bullet-\xi-n)\rangle|^2&=\|T_\xi^*f\|^2\\
			&=\|T_{-\xi}f\|^2\\
			&=\|T_\xi f\|^2, \  \text{for \ all} \ f\in L^2(\mathbb{R}^d).
			\end{align*}
		\end{exa}
		
		The following proposition gives a characterization of $K$- $\mathcal{W} \ \mathcal{H}$ frame in term of operator inequality.
		\begin{pro}
			Let $g\in L^2(\mathbb{R}^d)$. Then, $\{E_{C_m}T_{\mathbf{B}n}g\}_{m\in \mathbb{Z}, n \in \mathbb{Z}^d} $ is a $K$- $\mathcal{W} \ \mathcal{H}$
			frame with bounds $A$,$B$ if and only if $AKK^{*}\leq
			\mathcal{S}\leq BK^{*}K$, where  $\mathcal{S}$ is frame
			operator for $\{E_{C_m}T_{\mathbf{B}n}g\}_{m\in \mathbb{Z}, n \in \mathbb{Z}^d} $.
		\end{pro}
		\proof Since
		\begin{align*}
		&\{E_{C_m}T_{\mathbf{B}n}g\}_{m\in \mathbb{Z}, n \in \mathbb{Z}^d} \ \text{is
			a} \ K \  \mathcal{W}\text{-}\mathcal{H} \ \text{frame} \  \text{with frame operator} \ S\\
		&\Leftrightarrow A\|K^{*}f\|^2\leq\sum_{m\in \mathbb{Z}, n \in \mathbb{Z}^d}|\langle f,E_{C_m}T_{\mathbf{B}n}g\rangle|^2\leq B\|K f\|^2
		\text{for \ all } f
		\in L^2(\mathbb{R}^d)\\
		&\Leftrightarrow A\|K^{*}f\|^2\leq\langle
		\mathcal{S}f,f\rangle\leq B\|K f\|^2 \text{for \ all } f
		\in L^2(\mathbb{R}^d)\\
		&\Leftrightarrow A\langle K^{*}f,K^{*}f\rangle\leq\langle
		\mathcal{S}f,f\rangle\leq B\langle K f,K f\rangle \
		\text{for \ all }
		f\in L^2(\mathbb{R}^d)\\
		&\Leftrightarrow A\langle KK^{*}f,f\rangle\leq\langle
		\mathcal{S}f,f\rangle\leq B\langle K^{*}K f, f\rangle \
		\text{for \
			all }f\in L^2(\mathbb{R}^d)\\
		&\Leftrightarrow AKK^{*}\leq \mathcal{S}\leq
		BK^{*}K .
		\end{align*}
		The proposition is proved.
		\endproof
		The following theorem provides necessary and sufficient condition
		for a given system to be  $K$- Weyl-Heisenberg
		frame for $L^2(\mathbb{R}^d).$
		\begin{thm}\label{t3.4}
			Let $g \in L^2(\mathbb{R}^d)$, then, $\{E_{C_m}T_{\mathbf{B}n}g\}_{m\in \mathbb{Z}, n \in \mathbb{Z}^d} $ is a $K$- Weyl-Heisenberg frame
			for $L^2(\mathbb{R}^d)$ if and only if there exist a  bounded linear
			operator $L:l^2\bigoplus l^2\longrightarrow L^2(\mathbb{R}^d)$ such
			that
			\begin{align*}
			L(e_{mn})=E_{C_m}T_{Bn}g \ \text{and} \ Range(K)\subset
			Range(L),
			\end{align*}
			where $\{e_{mn}\}$ is an orthonormal basis for $l^2\bigoplus l^2.$
		\end{thm}
		
		\proof
		Suppose first  that $\{E_{C_m}T_{\mathbf{B}n}g\}_{m\in \mathbb{Z}, n \in \mathbb{Z}^d}$ is a  $K$- Weyl-Heisenberg frame
		for $L^2(\mathbb{R}^d)$. Then, we can find positive constant
		$\mathrm{a}_0, \mathrm{b}_0$ such that
		\begin{align*}
		\mathrm{a}_0\|K^{*}f\|^2\leq \sum_{m\in \mathbb{Z}, n \in \mathbb{Z}^d}|\langle
		f,E_{C_m}T_{Bn}g\rangle|^2 \leq \mathrm{b}_0\| f\|^2, \  \text{for \ all} \  f\in L^2(\mathbb{R}^d).
		\end{align*}
		Define $ \mathcal{S}:L^2(\mathbb{R}^d) \longrightarrow \ell^2
		\bigoplus \ell^2$ by
		\begin{align*}
		\mathcal{S}(f)=\sum_{m\in \mathbb{Z}, n \in \mathbb{Z}^d}\langle
		f,E_{C_m}T_{\mathbf{B}n}g\rangle e_{mn}
		\end{align*}
		Clearly, $\mathcal{S}$ is a well defined bounded linear
		operator on $L^2(\mathbb{R}^d)$.\\
		Now
		\begin{align*}
		\langle \mathcal{S}^{*}e_{mn},h\rangle &= \langle e_{mn},\mathcal{S}h\rangle\\
		&=\langle e_{mn}
		\sum_{m\in \mathbb{Z}, n \in \mathbb{Z}^d}\langle h,E_{C_m}T_{\mathbf{B}n}g\rangle
		e_{mn}\rangle\\
		&=\sum_{m\in \mathbb{Z}, n \in \mathbb{Z}^d}\overline{\langle h,E_{C_m}T_{\mathbf{B}n}g\rangle}\langle
		e_{mn},e_{mn}\rangle\\
		&=\overline{\langle h,E_{C_m}T_{\mathbf{B}n}g\rangle}\\
		&=\langle E_{C_m}T_{\mathbf{B}n}g,h\rangle , \ \text{for \ all} \ h \in L^2(\mathbb{R}^d).
		\end{align*}
		Therefore, $\mathcal{S}^{*}e_{mn}=E_{C_m}T_{\mathbf{B}n}g.$
		
		Also
		\begin{align*}
		\mathrm{a}_0\|K^{*}f\|^2\leq\sum _{m\in \mathbb{Z}, n \in \mathbb{Z}^d}|\langle
		f,\mathcal{S}^{*}e_{mn}\rangle|^2=\|\mathcal{S}f\|^2 , \text{ for \
			all} \ f\in L^2(\mathbb{R}^d).
		\end{align*}
		Thus, $\mathrm{a}_0KK^{*}\leq \mathcal{S}^{*}\mathcal{S}$
		, where $L=\mathcal{S}^{*}$. Also, by Theorem \ref{t2.3}, we have
		Range$(K)\subset$~Range$(L).$
		
		Conversely, assume that  $E_{C_m}T_{Bn}g=Le_{mn}$ where $L\in \mathcal{B}(l^2\bigoplus l^2, L^2(\mathbb{R}^d))$  and
		Range$(K)\subset$Range$(L)$.
		
		Since
		\begin{align*}
		\langle L^{*}f,h \rangle &=\langle L^{*}f,\sum_{m\in \mathbb{Z}, n \in \mathbb{Z}^d}a_{mn}e_{mn}\rangle\\
		&=\sum_{m\in \mathbb{Z}, n \in \mathbb{Z}^d} \overline{a_{mn}}\langle f,Le_{mn}\rangle\\
		&=\sum_{m\in \mathbb{Z}, n \in \mathbb{Z}^d} \overline{a_{mn}}\langle
		f,E_{C_m}T_{\mathbf{B}n}g\rangle\\
		&=\sum_{m\in \mathbb{Z}, n \in \mathbb{Z}^d} \overline{\langle
			h,e_{mn}\rangle}\langle
		f,E_{C_m}T_{\mathbf{B}n}g\rangle\\
		&=\sum_{m\in \mathbb{Z}, n \in \mathbb{Z}^d} \langle e_{mn},h\rangle\langle
		f,E_{C_m}T_{\mathbf{B}n}g\rangle\\
		&=\langle\sum_{m\in \mathbb{Z}, n \in \mathbb{Z}^d} \langle f,E_{C_m}T_{\mathbf{B}n}g\rangle
		e_{mn},h\rangle\\
		\end{align*}
		Therefore, $L^{*}f= \sum_{m\in \mathbb{Z}, n \in \mathbb{Z}^d} \langle f,
		E_{C_m}T_{\mathbf{B}n}g\rangle e_{mn}.$\\
		Now
		\begin{align}\label{e3.2}
		\sum_{m\in \mathbb{Z}, n \in \mathbb{Z}^d} |\langle f,E_{C_m}T_{\mathbf{B}n}g\rangle|^2 &= \sum_{m\in \mathbb{Z}, n \in \mathbb{Z}^d} |\langle f,Le_{mn}\rangle|^2\notag\\
		&=\sum_{m\in \mathbb{Z}, n \in \mathbb{Z}^d} |\langle L^{*}f,e_{mn}\rangle|^2 \notag\\
		&=\|L^{*}f\|^2 \notag\\
		&\leq \|L^{*}\|^2\|f\|^2, \text{for \ all} \ f \in L^2(\mathbb{R}^d)
		\end{align}

		By Theorem \ref{t2.3} and using fact that Range$(K)\subset$
		Range$(L).$, we can find a positive
		constant $A$  such  that $AKK^{*}\leq LL^{*}$.\\
		Therefore, for all $f\in L^2(\mathbb{R}^d),$ we have
		\begin{align}\label{e3.3}
		A\|K^{*}f\|^2 &\leq\|L^{*}f\|^2 \notag\\
		&=\sum_{m\in \mathbb{Z}, n \in \mathbb{Z}^d} |\langle f,E_{C_m}T_{Bn}g\rangle|^2.
		\end{align}
		By using \eqref{e3.2} and \eqref{e3.3}, we have
		\begin{align*}
		A\|K^{*}f\|^2\leq\sum_{m\in \mathbb{Z}, n \in \mathbb{Z}^d} |\langle
		f,E_{C_m}T_{Bn}g\rangle|^2\leq\|L^{*}\|^2\| f\|^2 , \ \text{for \
			all} \ f\in L^2(\mathbb{R}^d).
		\end{align*}
		Hence $\{E_{C_m}T_{Bn}g\}_{m\in \mathbb{Z}, n \in \mathbb{Z}^d} $ is a K-$\mathcal{W} \ \mathcal{H}$ frame  for
		$L^2(\mathbb{R}^d).$
		\endproof
		
		The following theorem provides sufficient conditions for the
		existence  of a $K$-$\mathcal{W} \ \mathcal{H}$ frame for $L^2(\mathbb{R}^d)$.
		\begin{thm}\label{t3.5}
			Let $K \in \mathcal{B}\left(L^2(\mathbb{R}^d)\right)$. Suppose
			that  $\{C_m\}_{m \in \mathbb{Z}}$ is a set of uniform density in
			$\mathbb{R}$ and  $g \in L^2(\mathbb{R}^d)$ is bounded with support
			$[-a, a]^d$, with
			\begin{align*}
			A\|K\|^2 \leq \sum_{n \in \mathbb{Z}^d}|g(t-\mathbf{B}n)|^2 \leq B\|K\|^{-2}  \ a.e. \quad(0 < A, B \in \mathbb{R})
			\end{align*}
			Then, $\{E_{C_m}T_{\mathbf{B}n}g\}_{m\in \mathbb{Z}, n \in \mathbb{Z}^d}$ is a
			$K$- $\mathcal{W} \ \mathcal{H}$ frame for $L^2(\mathbb{R}^d)$.
		\end{thm}
		\proof Since $\{C_m\}_{m \in \mathbb{Z}}$ is a set of uniform
		density in $\mathbb{R}^d$, there is an $a>0$ such that $\{E_{C_m}\}_{m
			\in \mathbb{Z}}$ is a frame for $L^2([-a,a]^d)$. Let   $A_1, \ B_1$ be
		a choice of  bounds for $\{E_{C_m}\}_{m \in \mathbb{Z}}$. Now for
		all $f \in L^2(\mathbb{R}^d)$, we have
		\begin{align}\label{e3.4}
		\sum_{m\in \mathbb{Z}, n \in \mathbb{Z}^d} |\langle f,E_{C_m}T_{\mathbf{B}n}g\rangle|^2 &= \sum_{m\in \mathbb{Z}, n \in \mathbb{Z}^d} |\langle
		f.T_{\mathbf{B}n}\overline{g},E_{C_m}\rangle|^2 \notag\\
		& \leq B_1\sum_{n \in \mathbb{Z}^d}\|f.T_{\mathbf{B}n}\overline{g}\|^2 \notag\\
		&=B_1\sum_{n \in \mathbb{Z}^d}\int_{\mathbb{R}^d}|f(t)|^2|g(t-\mathbf{B}n)|^2dt \notag\\
		&=B_1\int_{\mathbb{R}^d}|f(t)|^2\sum_{n \in \mathbb{Z}^d}|g(t-\mathbf{B}n)|^2dt \notag\\
		&\leq B_1B\|K\|^{-2}\int_{\mathbb{R}^d}|f(t)|^2dt \notag\\
		&\leq B_1B\|K\|^{-2}\|f\|^2.
		\end{align}
		Also
		\begin{align}\label{e3.5}
		A_1A\|K^{*} f\|^2 &\leq A_1A\|K^{*}\|^2\| f\|^2 \notag\\
		&=A_1A\|K\|^2\int_{\mathbb{R}^d}|f(t)|^2dt \notag\\
		&=A_1\int_{\mathbb{R}^d}|f(t)|^2\sum_{n \in \mathbb{Z}^d}|g(t-Bn)|^2dt \notag\\
		&=A_1\sum_{n \in \mathbb{Z}^d}\|f.T_{\mathbf{B}n}\overline{g}\|^2  \notag\\
		&=\sum_{n \in \mathbb{Z}^d}A_1\|f.T_{\mathbf{B}n}\overline{g}\|^2 \notag\\
		&\leq\sum_{n \in \mathbb{Z}^d}\sum_{m \in \mathbb{Z}}|\langle
		f.T_{\mathbf{B}n}\overline{g},E_{C_m}\rangle|^2 \notag\\
		&=\sum_{m\in \mathbb{Z}, n \in \mathbb{Z}^d} |\langle f,E_{C_m}T_{\mathbf{B}n}g\rangle|^2.
		\end{align}
		By using \eqref{e3.4} and \eqref{e3.5},  $\{E_{C_m}T_{\mathbf{B}n}g\}_{m\in \mathbb{Z}, n \in \mathbb{Z}^d}$ is a $K$-
		$\mathcal{W} \ \mathcal{H}$ frame for $L^2(\mathbb{R}^d)$ with bounds
		$A_1A$  and $B_1B \|K\|^{-2}$. This complete the proof.
		\endproof
		
		The following theorem gives necessary condition for the
		$K$- $\mathcal{W} \ \mathcal{H}$ frame.
		
		\begin{thm}
			Let $\{E_{C_m}T_{\mathbf{B}n}g\}_{m\in \mathbb{Z}, n \in \mathbb{Z}^d}$ is a
			$K$- $\mathcal{W} \ \mathcal{H}$ frame for
			$L^2(\mathbb{R}^d)$ with bounds $A$ and $B$. Then, there exists a
			positive number $\tau$ such that $\{E_{C_m}\}_{m  \in \mathbb{Z}}$
			is a frame for $L^2[0, \tau]^d$. Furthermore, if   $A_o, B_o$ are
			bounds of $\{E_{C_m}\}_{m  \in \mathbb{Z}}$, then
			\begin{align*}
			\sum_{n \in \mathbb{Z}^d}|g(t-\mathbf{B}n)|^2 \leq \frac{B}{A_o} \|K\|^2.
			\end{align*}
		\end{thm}
		\proof Let $\mathcal{J}= \{(C_m, Bn): m\in\mathbb{Z}, n \in \mathbb{Z}^d\}$. Then,
		there exist $r, \omega > 0$ such that for all $a, b \in \mathbb{R}$
		$([a ,a+r]^d \times [b, b+\omega]^d)\bigcap \mathcal{J}\neq \phi.$
		Therefore, $\{C_m\}_{m \in \mathbb{Z}}$ is a system of uniform
		density in $\mathbb{R}^d$. Thus,  there exists $\tau
		> 0$ such that $\{E_{C_m}\}_{m \in \mathbb{Z}}$ is a frame for
		$L^2[0, \tau]^d$. Let $A_o, B_o$ be a choice of bounds for
		$\{E_{C_m}\}_{m \in \mathbb{Z}}$. Then, for any interval $I= [c,
		c+\tau]^d$ and bounded function $f \in L^2(I)$, we have
		\begin{align*}
		\langle f,E_{C_m}T_{\mathbf{B}n}g\rangle &=\int_I f(t)\overline{E_{C_m}T_{Bn}g(t)} dt\\
		&=\int_I f(t)\overline{e^{2\pi iC_m.t}g(t-\mathbf{B}n)}dt\\
		&=\int_I f(t)\overline{g}(t-\mathbf{B}n) \overline{e^{2\pi iC_m.t}}dt\\
		&=\int_I (f.T_{\mathbf{B}n}\overline{g})(t) \overline{E_{C_m}(t)}dt\\
		&= \langle f.T_{\mathbf{B}n}\overline{g}, E_{C_m} \rangle.
		\end{align*}
		Therefore
		\begin{align*}
			\sum_{n \in \mathbb{Z}^d}\|f.T_{\mathbf{B}n}\overline{g}\|^2  &\leq\sum_{m\in \mathbb{Z}, n \in \mathbb{Z}^d}|\langle f.T_{\mathbf{B}n}\overline{g}, E_{C_m}\rangle|^2\\
			&=\sum_{m\in \mathbb{Z}, n \in \mathbb{Z}^d} |\langle
			f,E_{C_m}T_{\mathbf{B}n}g\rangle|^2\\
			&\leq B\|K f\|^2\\
			&\leq B \|K\|^2\|f\|^2.
		\end{align*}
		Thus
		\begin{align*}
		A_o\int_I|f(t)|^2\sum_{n \in \mathbb{Z}^d}|g(t-\mathbf{B}n)|^2dt\leq B
		\|K\|^2\int_I|f(t)|^2 \text{ for \ all } f \in
		L^2(\mathbb{R}^d).\\
		\end{align*}
		Hence $\sum_{n \in \mathbb{Z}^d}|g(t-\mathbf{B}n)|^2\leq \frac{B}{A_o}
		\|K\|^2.$
		\endproof
		It is interesting to know whether a given system in the structure of
		algebra and topology, is invariant under certain operator. Next two
		results in this section discuss the image of a given standard
		$\mathcal{W} \ \mathcal{H}$ frame for $L^2(\mathbb{R}^d)$. First
		result show that for a given $K \in
		\mathcal{B}(L^2(\mathbb{R}^d))$, the image of a standard $\mathcal{W} \ \mathcal{H}$ frame for $L^2(\mathbb{R}^d)$ is a
		$\mathcal{W} \ \mathcal{H}$ frame  for $L^2(\mathbb{R}^d)$.
		Second result reflects the invariance of frame properties of
		$K$- $\mathcal{W} \ \mathcal{H}$ frame under a linear
		homeomorphism, together with a relation with best bounds of the
		original frame.
		
		\begin{thm}
			Let $g\in L^2(\mathbb{R}^d)$ and let $\{E_{C_m}T_{\mathbf{B}n}g\}_{m\in \mathbb{Z}, n \in \mathbb{Z}^d}$ be a standard $\mathcal{W} \ \mathcal{H}$-frame
			frame for $L^2(\mathbb{R}^d)$. Then, $\mathcal{F}_{K} \equiv
			\{K (E_{C_m}T_{\mathbf{B}n}g)\}_{m\in \mathbb{Z}, n \in \mathbb{Z}^d}$ is a
			$K$- $\mathcal{W} \ \mathcal{H}$ frame for
			$L^2(\mathbb{R}^d)$.
		\end{thm}
		\proof Since $\{E_{C_m}T_{\mathbf{B}n}g\}_{m\in \mathbb{Z}, n \in \mathbb{Z}^d}$ is a
		ordinary $\mathcal{W} \ \mathcal{H}$-frame for
		$L^2(\mathbb{R}^d)$ Therefore there exists positive constants $A$,$B$
		such that
		\begin{align*}
		A\|f\|^2\leq \sum_{m\in \mathbb{Z}, n \in \mathbb{Z}^d}|\langle
		f,E_{C_m}T_{\mathbf{B}n}g\rangle|^2 \leq B\|f\|^2, \  \text{for   all} \
		f\in L^2(\mathbb{R}^d).
		\end{align*}
		Now
		\begin{align}\label{e3.6}
		\sum_{m\in \mathbb{Z}, n \in \mathbb{Z}^d}|\langle
		f,K(E_{C_m}T_{\mathbf{B}n}g)\rangle|^2 &=\sum_{m\in \mathbb{Z}, n \in \mathbb{Z}^d}|\langle K^{*}
		f,E_{C_m}T_{\mathbf{B}n}g\rangle|^2\notag\\
		&\leq B\|K^{*} f\|^2\notag\\
		&\leq B\|K\|^2\| f\|^2 \  \text{for   all} \ f\in
		L^2(\mathbb{R}^d).
		\end{align}
		Also
		\begin{align}\label{e3.7}
		A\|K^{*} f\|^2 & \leq\sum_{m\in \mathbb{Z}, n \in \mathbb{Z}^d}|\langle
		K^{*} f,E_{C_m}T_{\mathbf{B}n}g\rangle|^2\notag\\
		&=\sum_{m\in \mathbb{Z}, n \in \mathbb{Z}^d}|\langle
		f,K(E_{C_m}T_{\mathbf{B}n}g)\rangle|^2 \  \text{for   all} \ f\in L^2(\mathbb{R}^d).
		\end{align}
		Use inequality \eqref{e3.6} and \eqref{e3.7} we have
		\begin{align*}
		A\|K^{*} f\|^2\leq\sum_{m\in \mathbb{Z}, n \in \mathbb{Z}^d}|\langle
		f,K(E_{C_m}T_{\mathbf{B}n}g)\rangle|^2\leq B\|K\|^2\|f\|^2 \
		\text{for   all} \ f\in L^2(\mathbb{R}^d).
		\end{align*}
		Hence, $\mathcal{F}_{K}$ is a $K$-
		$\mathcal{W} \ \mathcal{H}$ frame for $L^2(\mathbb{R}^d)$ This
		completes the proof.
		\endproof
		
		The following theorem show that $K$-$\mathcal{W} \ \mathcal{H}$ frame for $L^2(\mathbb{R}^d)$ are
		invariant under a linear homeomorphism, provided both $K$ and
		its conjugate commutes with the given homeomorphism. A relation
		between the best bounds of a given $K$-
		$\mathcal{W} \ \mathcal{H}$ frame and best bounds of
		$K$- $\mathcal{W} \ \mathcal{H}$ frame  obtained by
		the action of linear homeomorphism is given.
		
		\begin{thm}
			Let $g\in L^2(\mathbb{R}^d)$. Suppose that $\{E_{C_m}T_{\mathbf{B}n}g\}_{m\in \mathbb{Z}, n \in \mathbb{Z}^d}$ is a $K$- $\mathcal{W} \ \mathcal{H}$
			frame for $L^2(\mathbb{R}^d)$ with best bounds $A_1$ and $B_1$. If
			$U:L^2(\mathbb{R}^d) \longrightarrow L^2(\mathbb{R}^d)$  is a linear
			homeomorphism such that $U$ commutes with both $K$ and
			$K^{*}$, then $\{U(E_{C_m}T_{\mathbf{B}n}g)\}_{m\in \mathbb{Z}, n \in \mathbb{Z}^d}$ is
			a $K$-frame for $L^2(\mathbb{R}^d)$ and its best bounds $A_2$,
			$B_2$ satisfy the inequalities
			\begin{align}\label{e3.8}
			A_1\|U\|^{-2}\leq A_2\leq A_1\|U^{-1}\|^2  \ ; \ B_1\|U\|^{-2}\leq B_2\leq B_1\|U\|^2.
			\end{align}
		\end{thm}
		\proof Since $B_1$ is an upper bound for $\{E_{C_m}T_{\mathbf{B}n}g\}_{m\in \mathbb{Z}, n \in \mathbb{Z}^d}$. Therefore, for all $f\in~L^2(\mathbb{R}^d)$, we have
		\begin{align}\label{e3.9}
		\sum_{m\in \mathbb{Z}, n \in \mathbb{Z}^d}|\langle f,U(E_{C_m}T_{Bn}g)\rangle|^2
		&=\sum_{m\in \mathbb{Z}, n \in \mathbb{Z}^d}|\langle
		U^{*}f,E_{C_m}T_{Bn}g\rangle|^2 \notag\\
		&\leq B_1\|K U^{*}f\|^2 \notag\\
		&\leq B_1\|U^{*}\|^2\|K f\|^2.
		\end{align}
		Also, by using the fact that $A_1$ is one of the choice for lower
		bound for $\{E_{C_m}T_{\mathbf{B}n}g\}_{m,n \in \mathbb{Z}}$, we have
		\begin{align}\label{e3.10}
		\|K^{*}f\|^2 &=\|K^{*}(UU^{-1})f\|^2 \notag\\
		&=\|UK^{*}(U^{-1}f)\|^2 \notag\\
		&\leq\|U\|^2\|K^{*}(U^{-1}f)\|^2 \notag\\
		&\leq\frac{\|U\|^2}{A_1}\sum_{m\in \mathbb{Z}, n \in \mathbb{Z}^d}|\langle
		U^{-1}f,E_{C_m}T_{\mathbf{B}n}g\rangle|^2  \notag\\
		&=\frac{\|U\|^2}{A_1}\sum_{m\in \mathbb{Z}, n \in \mathbb{Z}^d}|\langle
		UU^{-1}f,U(E_{C_m}T_{\mathbf{B}n}g)\rangle|^2\notag\\
		&=\frac{\|U\|^2}{A_1}\sum_{m\in \mathbb{Z}, n \in \mathbb{Z}^d}|\langle
		f,U(E_{C_m}T_{\mathbf{B}n}g)\rangle|^2.
		\end{align}
		By using \eqref{e3.9} and \eqref{e3.10}, we obtain
		\begin{align*}
		&A_1\|U\|^{-2}\|K^{*}f\|^2\leq\sum_{m\in \mathbb{Z}, n \in \mathbb{Z}^d}|\langle f,U(E_{C_m}T_{\mathbf{B}n}g)\rangle|^2\leq B_1
		\|U^*\|\|K f\|^2,\ \text{\ for \ all} \ f\in L^2(\mathbb{R}^d).
		\end{align*}
		Hence $\{U(E_{C_m}T_{\mathbf{B}n}g)\}_{m\in \mathbb{Z}, n \in \mathbb{Z}^d}$ is a
		$K$- $\mathcal{W} \ \mathcal{H}$ frame for
		$L^2(\mathbb{R}^d)$ with one of the choice of a pair $(A_1\|U\|^{-2},
		B_1\|U\|^2)$ as its bounds. Since $A_2$ and $B_2$ are best bounds
		for $\{U(E_{C_m}T_{\mathbf{B}n}g)\}_{m\in \mathbb{Z}, n \in \mathbb{Z}^d}$, we have
		\begin{align}\label{e3.11}
		A_1\|U\|^{-2}\leq A_2 , B_2\leq B_1\|U\|^2.
		\end{align}
		Since $\{U(E_{C_m}T_{\mathbf{B}n}g)\}_{m\in \mathbb{Z}, n \in \mathbb{Z}^d}$ is a
		$K$- $\mathcal{W} \ \mathcal{H}$ frame for
		$L^2(\mathbb{R}^d)$ with $(A_2, B_2)$ as one of the choice of bounds.
		So, for all $f \in L^2(\mathbb{R}^d),$ we have
		\begin{align}\label{e3.12}
		A_2\|K^{*}f\|^2\leq\sum_{m\in \mathbb{Z}, n \in \mathbb{Z}^d}|\langle
		f,U(E_{C_m}T_{\mathbf{B}n}g)\rangle|^2\leq B_2\|K f\|^2.
		\end{align}
		Now
		\begin{align}\label{e3.13}
		\|K^{*}f\|^2 &=\|U^{-1}UK{*}f\|^2 \notag\\
		&=\|U^{-1}K^{*}Uf\|^2 \notag\\
		&\leq\|U^{-1}\|^2\|K^{*}Uf\|^2.
		\end{align}
		By using \eqref{e3.12} and \eqref{e3.13}, we have
		\begin{align}\label{e3.14}
		A_2\|U^{-1}\|^{-2}\|K^{*}f\|^2  &\leq A_2\|K^{*}Uf\|^2 \notag \\
		&\leq \sum_{m\in \mathbb{Z}, n \in \mathbb{Z}^d}|\langle
		Uf,U(E_{C_m}T_{\mathbf{B}n}g)\rangle|^2 (=\sum_{m\in \mathbb{Z}, n \in \mathbb{Z}^d}|\langle f,E_{C_m}T_{\mathbf{B}n}g\rangle|^2) \notag\\
		&\leq B_2\|K Uf\|^2 \notag\\
		&\leq B_2\|U\|^2\|K f\|^2 ,\text{ \ for \ all} \ f\in
		L^2(\mathbb{R}^d).
		\end{align}
		Now $A_1$ and $B_1$ are the best bounds for $\{E_{C_m}T_{\mathbf{B}n}g\}_{m\in \mathbb{Z}, n \in \mathbb{Z}^d}$. Therefore, by using \eqref{e3.14}, we have
		\begin{align}\label{e3.15}
		A_2\|U^{-1}\|^{-2}\leq A_1 , B_1\leq B_2\|U\|^2.
		\end{align}
		The inequalities in \eqref{e3.8} is now obtained by using \eqref{e3.11} and \eqref{e3.15}.
		\endproof

	\section{Conclusion}
		We studied the notion of $K$-Weyl-Heisenberg frames (or $K$-$\mathcal{W} \ \mathcal{H}$ frames ) in $L^2(\mathbb{R}^d)$. Necessary and
		sufficient conditions for a certain system in $L^2(\mathbb{R}^d)$ to $K$-$\mathcal{W} \ \mathcal{H}$  frames for $L^2(\mathbb{R}^d)$ are obtained.
		We showed that $K$-$\mathcal{W} \ \mathcal{H}$ frame for $L^2(\mathbb{R}^d)$ are
		invariant under a linear homeomorphism, provided both $K$ and
		its conjugate commutes with the given homeomorphism. A relation
		between the best bounds of a given $K$-
		$\mathcal{W} \ \mathcal{H}$ frame and best bounds of
		$K$- $\mathcal{W} \ \mathcal{H}$ frame  obtained by
		the action of linear homeomorphism is given.
	
\end{document}